\documentclass[12pt]{article}
\usepackage{amsxtra}
\usepackage{amssymb, amsmath}
\usepackage{hyperref}
\usepackage {graphicx}

\begin{document}
\begin{center}
\textbf{\LARGE{Equivalent  Versions of ``\textit{Khajuraho}'' and ``\textit{Lo-Shu}'' Magic Squares and  the day $\bf{1^{st}}$ October 2010 (01.10.2010)}}
\end{center}

\bigskip
\begin{center}
\textbf{\large{Inder Jeet Taneja}}\\
Departamento de Matem\'{a}tica\\
Universidade Federal de Santa Catarina\\
88.040-900 Florian\'{o}polis, SC, Brazil.\\
\textit{e-mail: ijtaneja@gmail.com\\
http://www.mtm.ufsc.br/$\sim$taneja}
\end{center}

\begin{abstract}
\textit{In this short note we shall give connection between the most perfect ``Khajuraho''  magic square of order $4\times 4$ discovered in 10th century and the ``Lo-Shu'' magic square of order $3 \times 3$ with the day October 1, 2010, i.e., 01.10.2010. The day has only three digits 0, 1 and 2. Here we have given an equivalent version of Khajuraho magic square using only three digits 0, 1 and 2. If we write the above date date in two parts, \textbf{\textit{0110 2010}} , interestingly, the sum of new magic square is the first part, i.e.  0110, and the numbers appearing in the magic square are from the second part. An equivalent version of ``Lo-Shu'' magic square of order $3 \times 3$ is also given }
\end{abstract}

\section{History}

The study of magic squares is very old in history.  The earliest known magic square is Chinese, recorded around 2800 B.C. described as "Lo-Shu" magic square.  (not sure, because many places we can see written as 2200 B.C. or 2500 B.C., etc.). It is a typical $3 \times 3$ magic square, where the numbers were represented by patterns not numerals. See below

\begin{center}
\includegraphics[bb=0mm 0mm 208mm 296mm, width=65.3mm, height=43.3mm, viewport=3mm 4mm 205mm 292mm]{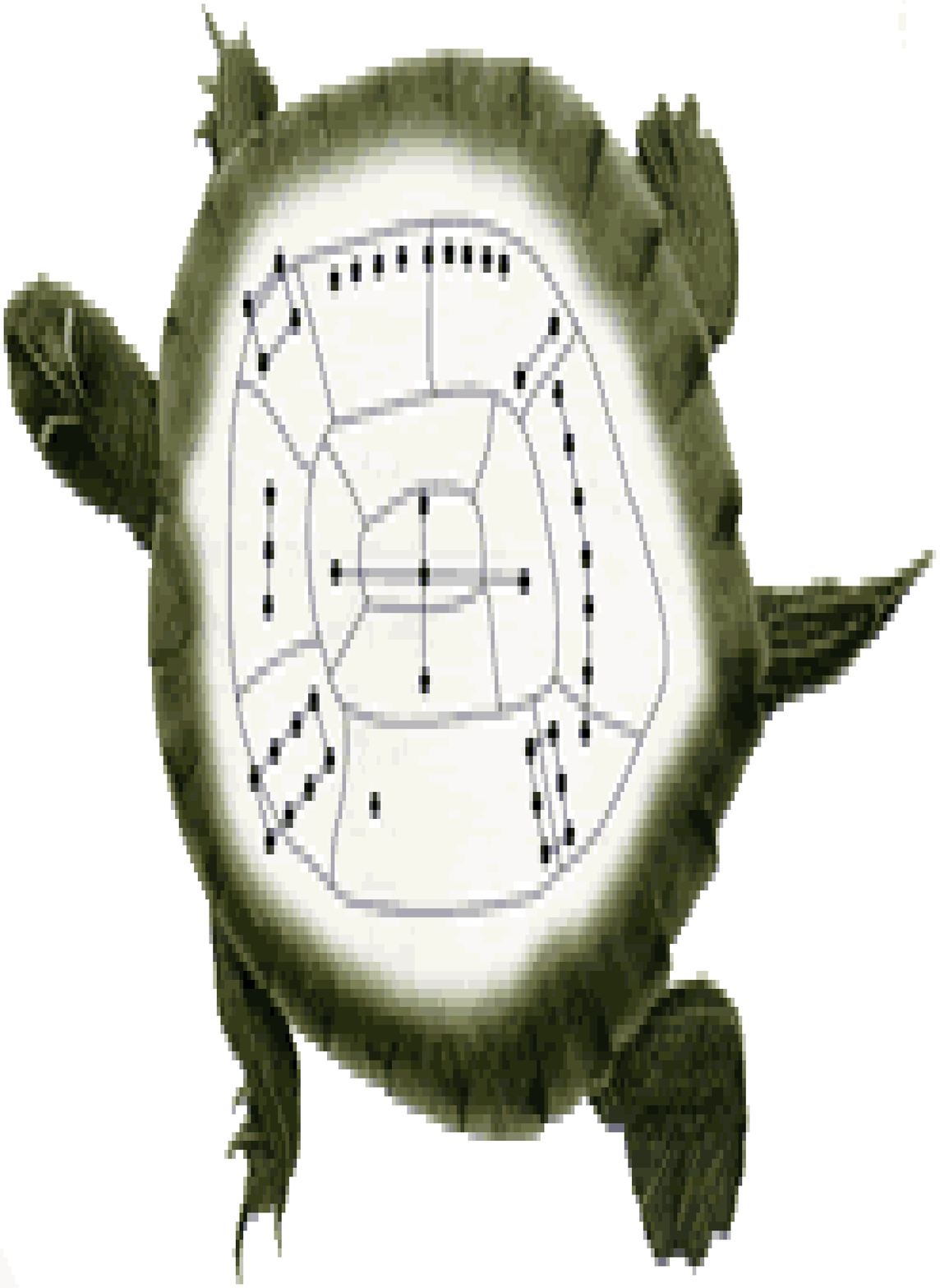}
\end{center}

(link: http://illuminations.nctm.org/LessonDetail.aspx?id=L263)
 -- accessed on 01.11.2010.

\bigskip
More precisely it is like this

\begin{center}
\includegraphics[bb=0mm 0mm 208mm 296mm, width=58.1mm, height=46.1mm, viewport=3mm 4mm 205mm 292mm]{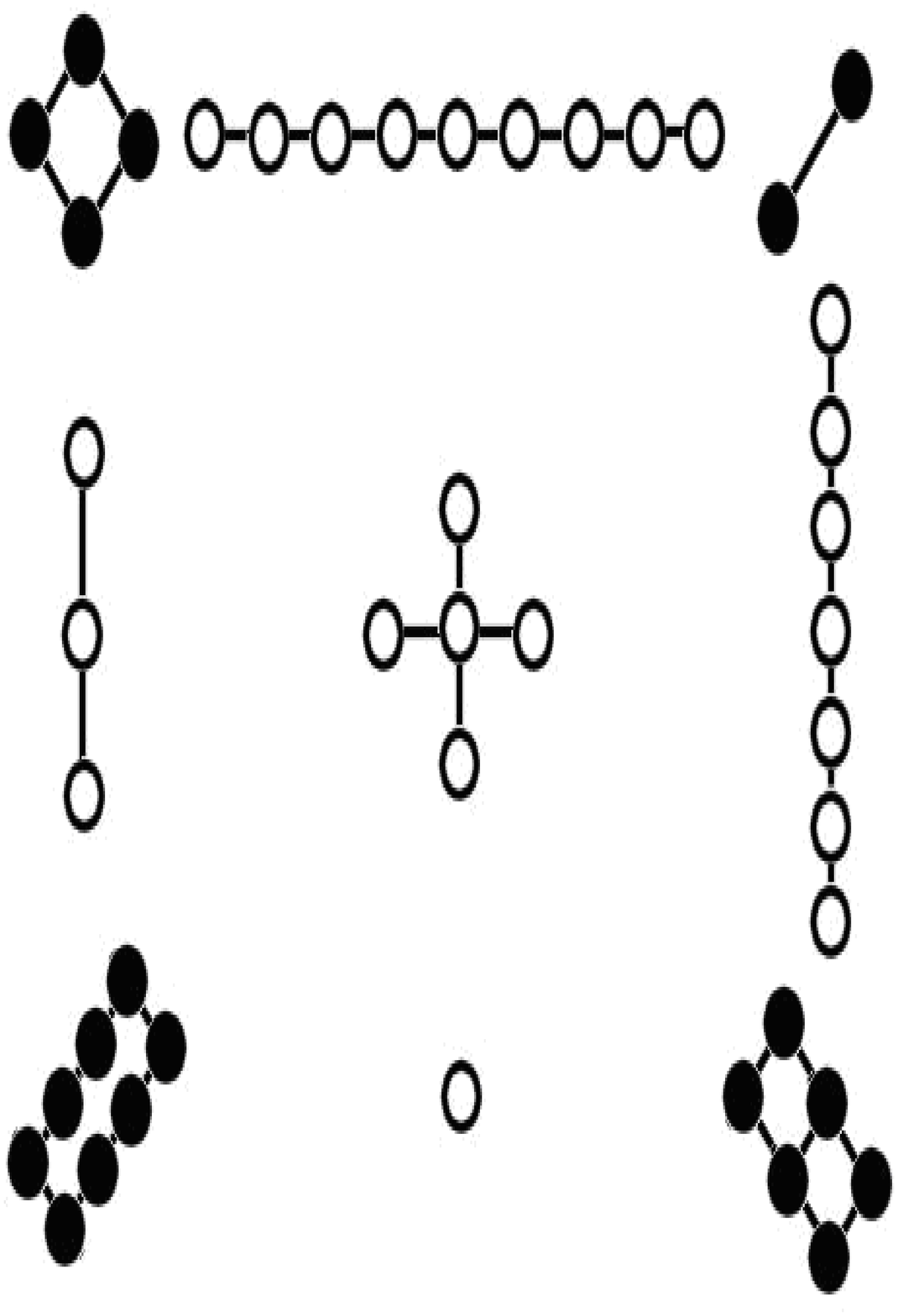}
\end{center}

Numerical transcription of above ``Lo-Shu'' magic square is

\begin{center}
\includegraphics[bb=0mm 0mm 208mm 296mm, width=22.2mm, height=20.7mm, viewport=3mm 4mm 205mm 292mm]{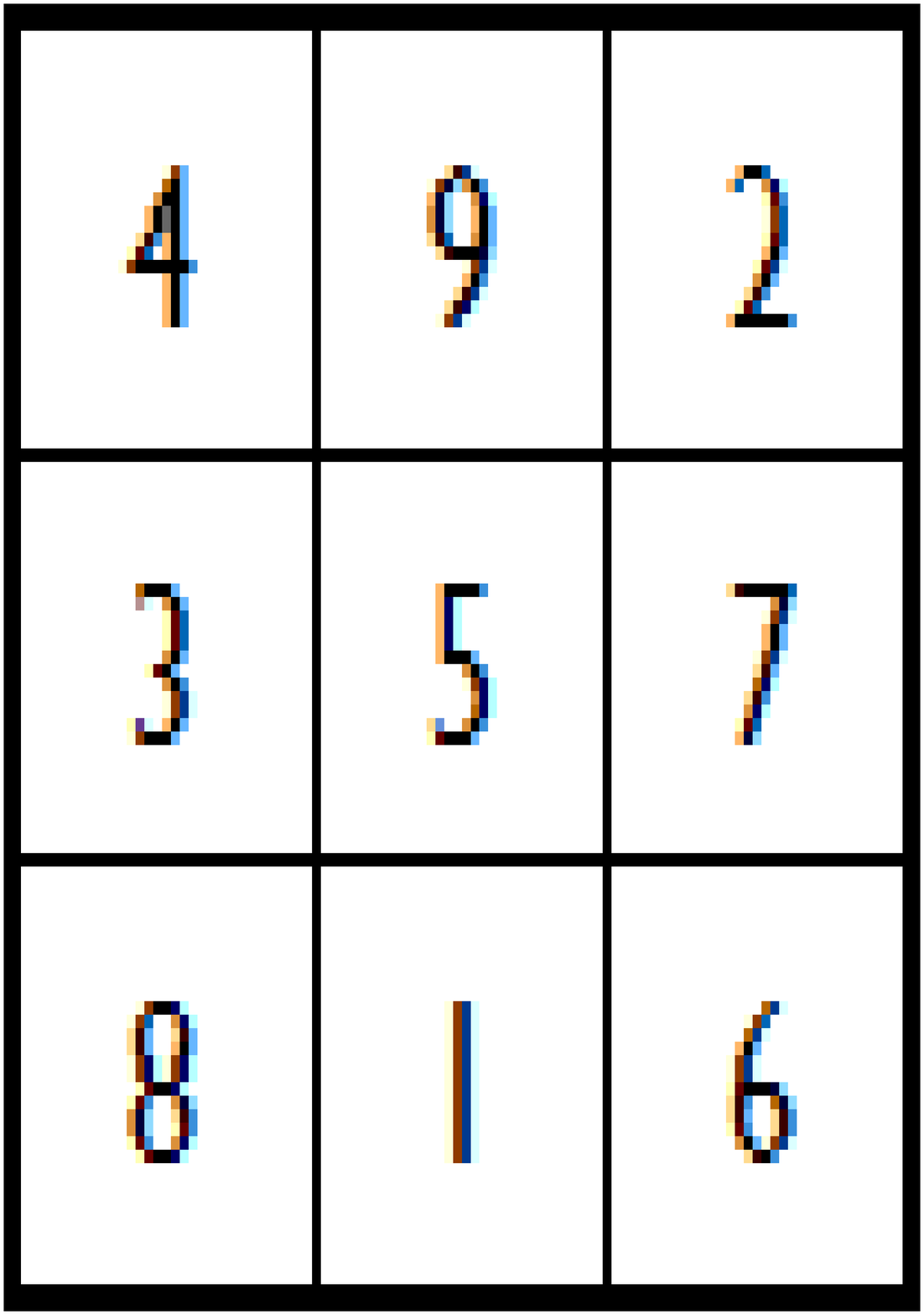}
\end{center}

Magic squares most likely traveled from China to India, then to the Arab countries. From the Arab countries, magic squares journeyed to Europe, then to Japan. Magic squares in India served multiple purposes other than the dissemination of mathematical knowledge. For example, Varahamihira used a fourth-order magic square to specify recipes for making perfumes in his book on seeing into the future,~\textit{Brhatsamhita}~(ca. 550 A.D.). The oldest dated third-order magic square in India appeared in Vrnda's medical work~\textit{Siddhayoga}~(ca. 900 A.D.), as a means to ease childbirth.~The fourth-order most perfect magic square found in $10^{th}$ century (945 AD) in Khajuraho Jain Mandir or Parshvanath Jain Tample.  See below plate appeared in Parshvanath Jain Tample in Khajuraho.

\begin{center}
\includegraphics[bb=0mm 0mm 208mm 296mm, width=48.4mm, height=73.7mm, viewport=3mm 4mm 205mm 292mm]{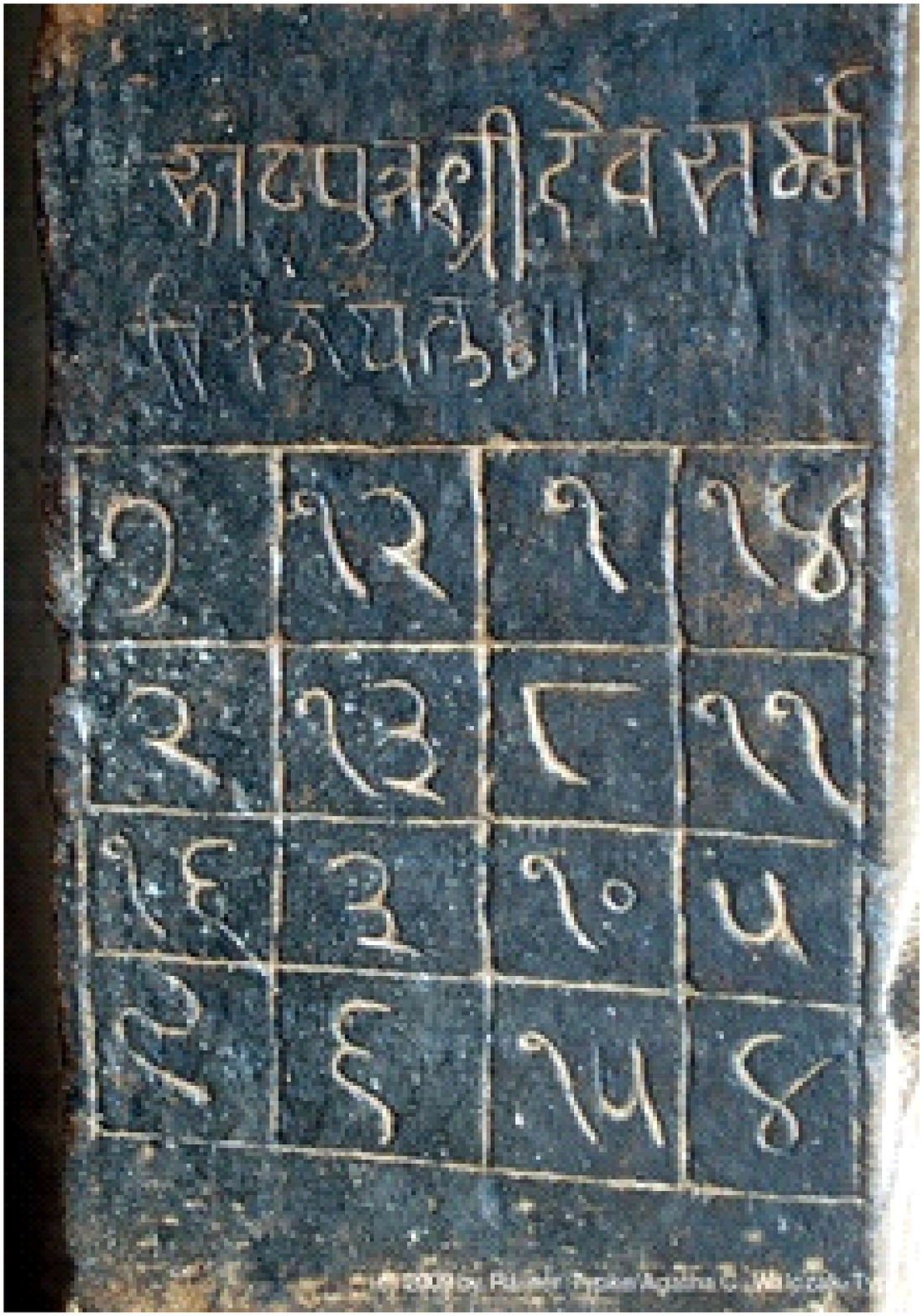}
\end{center}
(link: http://en.wikipedia.org/wiki/Jain\_temples\_of\_Khajuraho)
-- accessed on 01.11.2010.

\bigskip

It is one of the most perfect magic square of order 4x4. It is generally famous as \textit{Khajuraho magic square}. English transcription of above numbers is as follows:

\begin{center}
\includegraphics[bb=0mm 0mm 208mm 296mm, width=28.8mm, height=28.2mm, viewport=3mm 4mm 205mm 292mm]{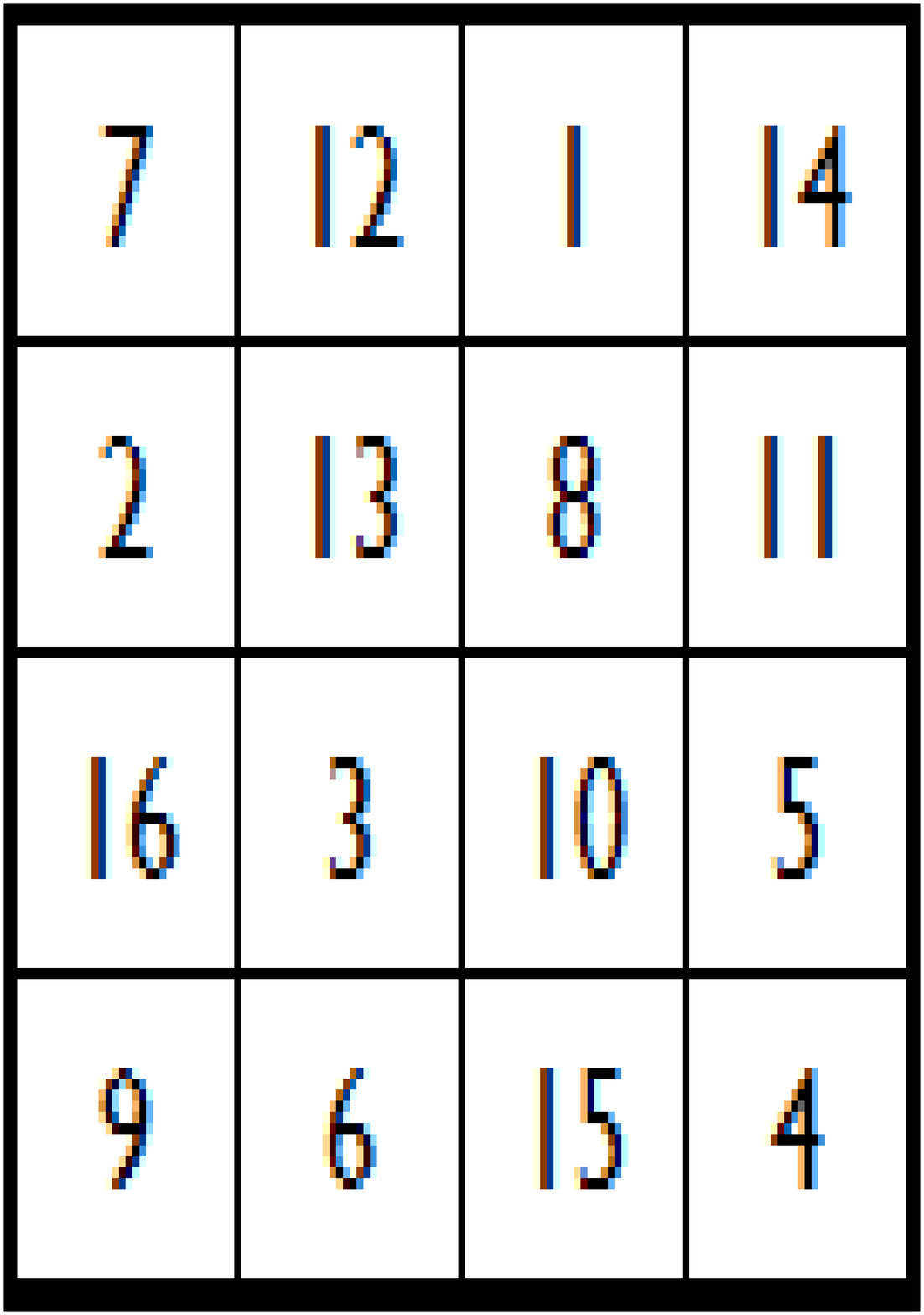}
\end{center}

It is pan diagonal magic square of sum 34. Each sub-square sum of 2x2 is also 34. Sometimes it is referred to as the \textit{``Chautisa Yantra''.}

\bigskip
We are now in $21^{st}$ century in a digital era, where everything is digitalized. The aim of this work is to produce an \textbf{\textit{equivalent upside down version}} of above magic square using only three digits appearing in a day 01.10.2010. In order to do so, we have used the numbers in the digital form: \textit{}

\begin{center}
\textit{\includegraphics[bb=0mm 0mm 208mm 296mm, width=18.6mm, height=5.6mm, viewport=3mm 4mm 205mm 292mm]{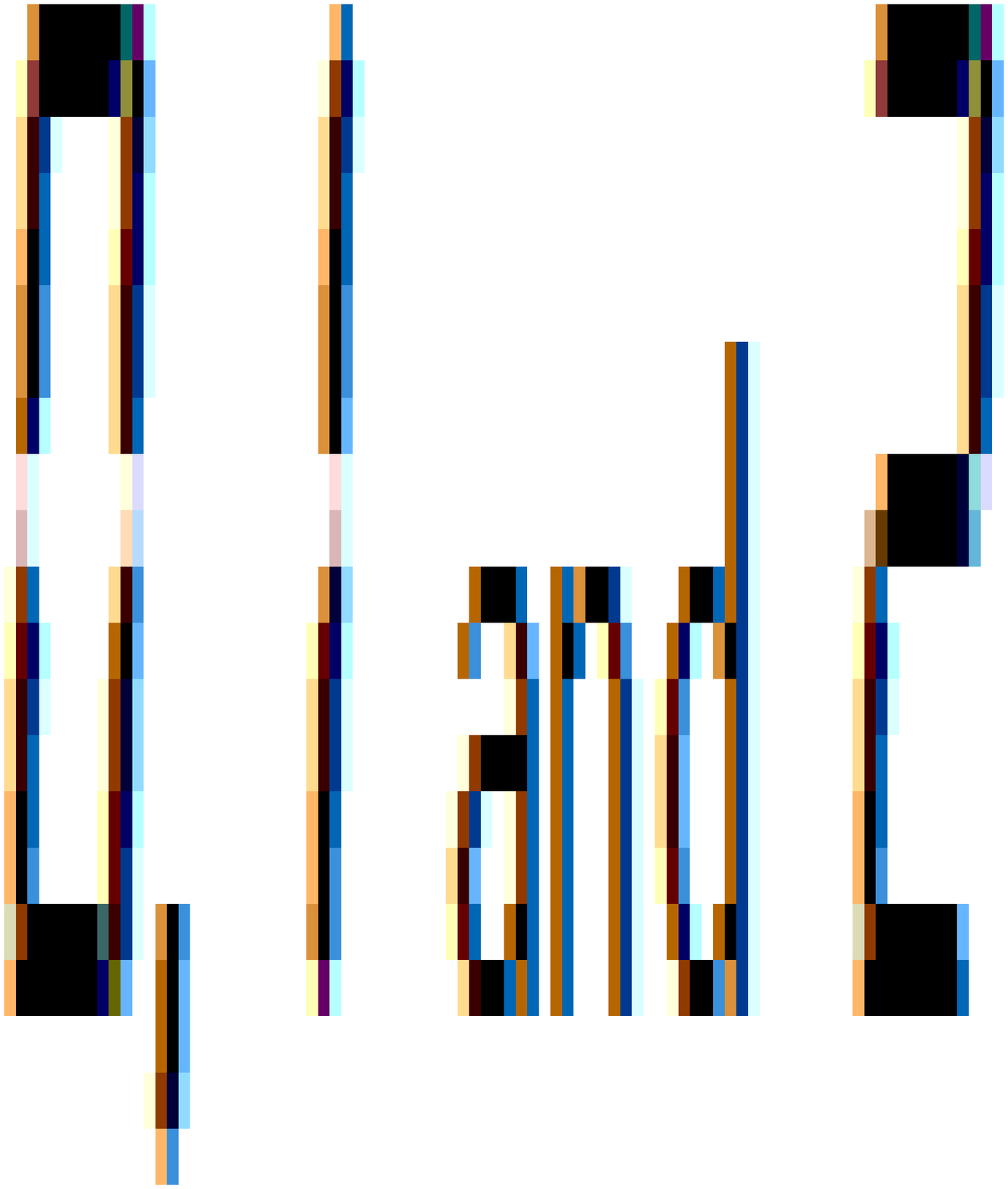}}
\end{center}

These digits generally appear in watches, elevators, etc.  We observe that the above three digits are rotatable to $180^{o}$, and remains the same.

\bigskip

Just to remember

\begin{itemize}
\item  [(i)]\textbf{A magic square} is a collection of numbers put as a square matrix, where the sum of elements of each row, sum of elements of each column and the sum of elements of each two principal diagonals are de same.

\item [(ii)]\textbf{ Upside down magic square }is a magic square, if we rotate it to $180^{o}$ (degrees) it remains again the magic square.

\item  [(iii)]\textbf{Mirror looking}, i.e., if we put it in front of mirror or see from the other side of the glass, or see on the other side of the paper, it always remains the magic square.
\end{itemize}

When the magic square is upside down and mirror looking, we call it \textbf{universal magic square}.

\section {Equivalent Upside Down Magic Squares}

\bigskip
In this section, we shall give \textbf{\textit{equivalent upside down versions}} of above two historical magic squares using only the digits 0, 1 and 2 in digital forms.

\subsection{Symmetric and Upside Down Version of  ``Khajuraho'' Magic Square}

\bigskip
Following Euler's idea of 1782 \cite{eul} of Latin squares, let us consider the following two \textbf{\textit{mutually orthogonal diagonalize Latin squares}} of order $4\times 4$.

\begin{center}
\includegraphics[bb=0mm 0mm 208mm 296mm, width=37.0mm, height=22.9mm, viewport=3mm 4mm 205mm 292mm]{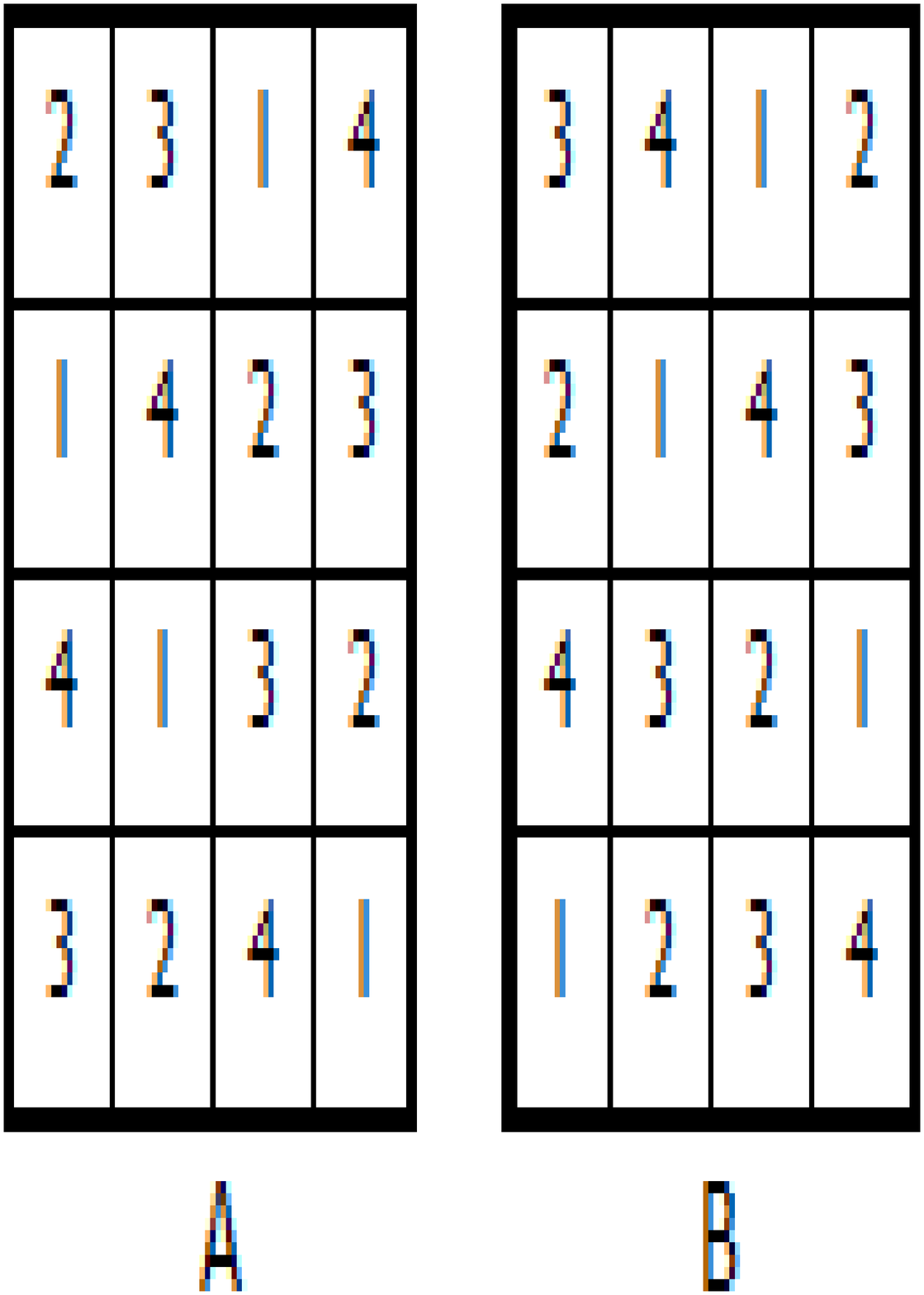}
\end{center}

  Let us apply an operation $4\times (A-1)+B$ in the above Latin squares, we get the following well know \textit{Khajuraho magic square} of order 4x4 appearing in the above plate:

\begin{center}
\includegraphics[bb=0mm 0mm 208mm 296mm, width=28.8mm, height=25.2mm, viewport=3mm 4mm 205mm 292mm]{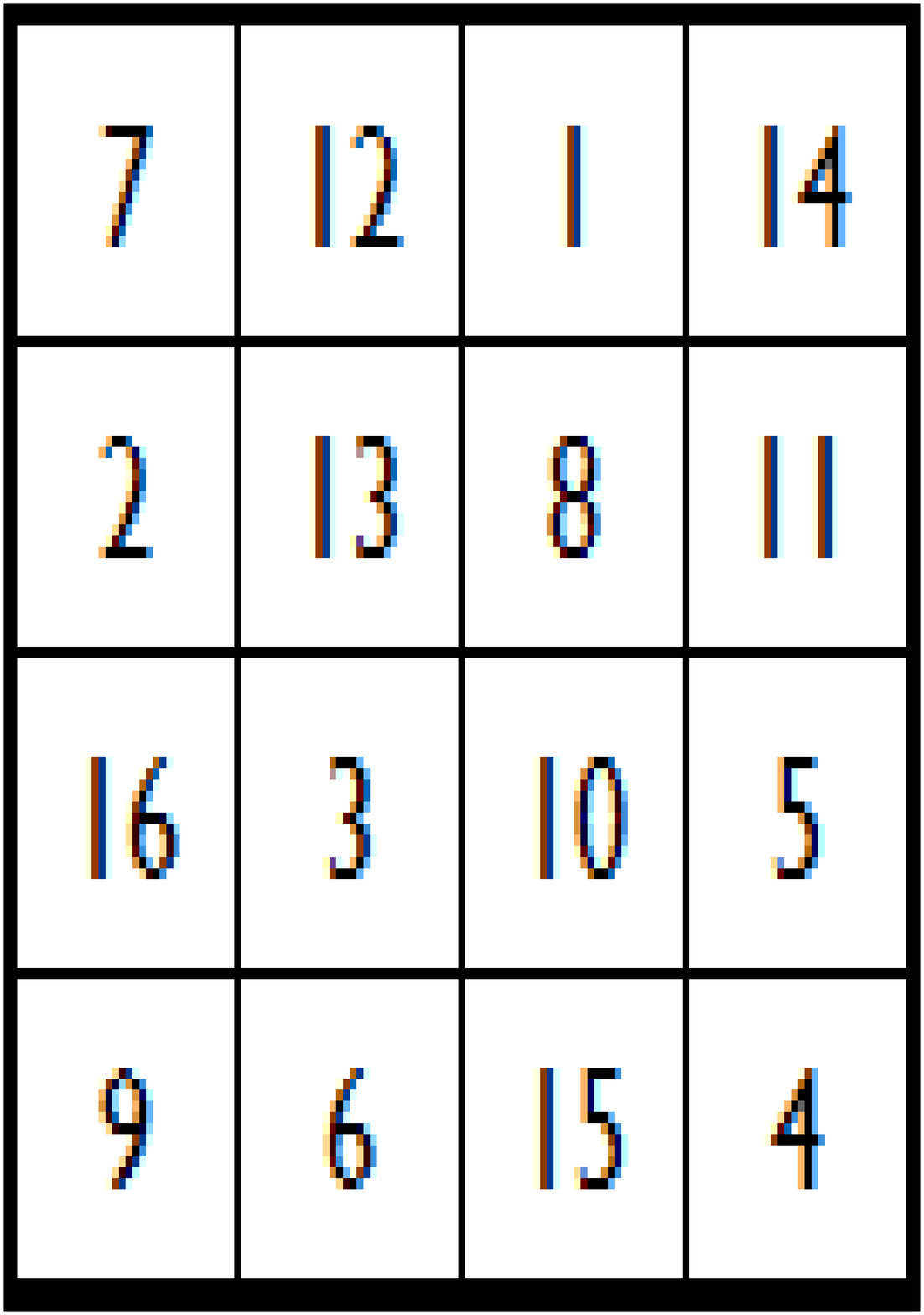}
\end{center}

As we wrote before the above magic square is pan diagonal. Also each sub-square of order 2x2 sums to 34.
\bigskip

Again if we apply an operation $10\times A+B$  over the above two Latin squares of order 4, then we get the following equivalent version of the \textit{Khajuraho magic square} with  sum as 110. See below

\begin{center}
\includegraphics[bb=0mm 0mm 208mm 296mm, width=41.7mm, height=38.7mm, viewport=3mm 4mm 205mm 292mm]{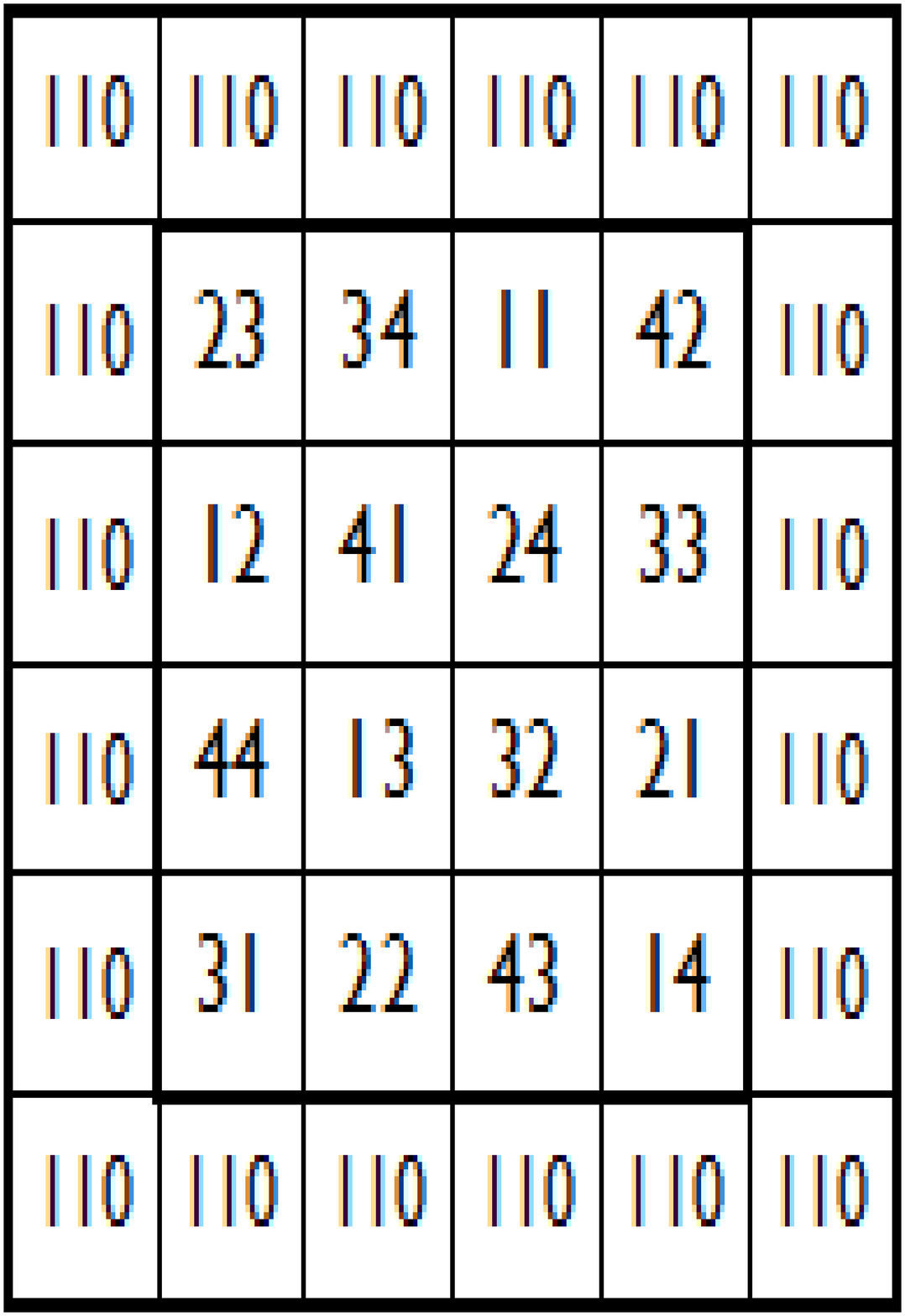}
\end{center}

The above magic square is not upside down but is \textbf{symmetric}, i.e, if we have \textit{ab, }then \textit{ba} is also there. Sometimes, we call it base 10 equivalent version.

\bigskip
Making some adjustments and writing the numbers in the digital form, here below is an \textbf{upside down equivalent version} of \textit{Khajuraho magic square}:

\begin{center}
\includegraphics[bb=0mm 0mm 208mm 296mm, width=85.9mm, height=52.7mm, viewport=3mm 4mm 205mm 292mm]{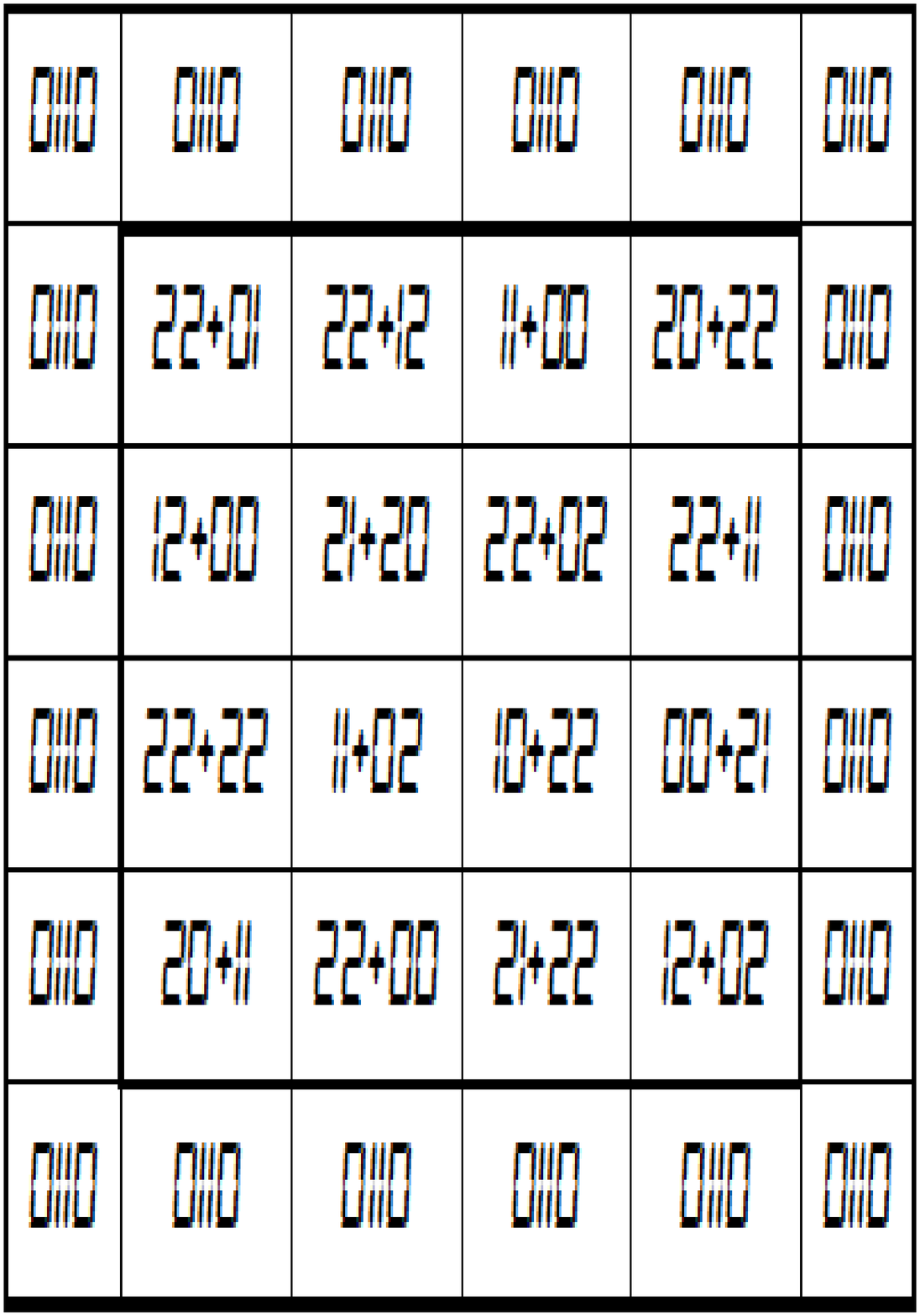}
\end{center}

The above magic square has only three digits 0, 1 and 2. If we gave a rotation of $180^{o}$ degrees it remains the same. Interestingly, the sum 110 writing as 0110 becomes $180^{o}$ degrees rotatable. We can consider it \textit{universal magic square}, since looking from the mirror we again get a magic square, where 2 becomes as 5. Naturally in this case, the sum is not the same. We leave it to readers to check it.

\bigskip
\noindent \textbf{$\bullet$ \large{The day 01.11.10}}

\bigskip
Instead, considering 2010, if we consider a year as 10, then the day we write as 01.11.10. It becomes a number 011110  (naturally we are putting 0 in the front to make it symmetric). The palindromic version of above magic square in four algorism is given by

\begin{center}
\includegraphics[bb=0mm 0mm 208mm 296mm, width=82.1mm, height=46.6mm, viewport=3mm 4mm 205mm 292mm]{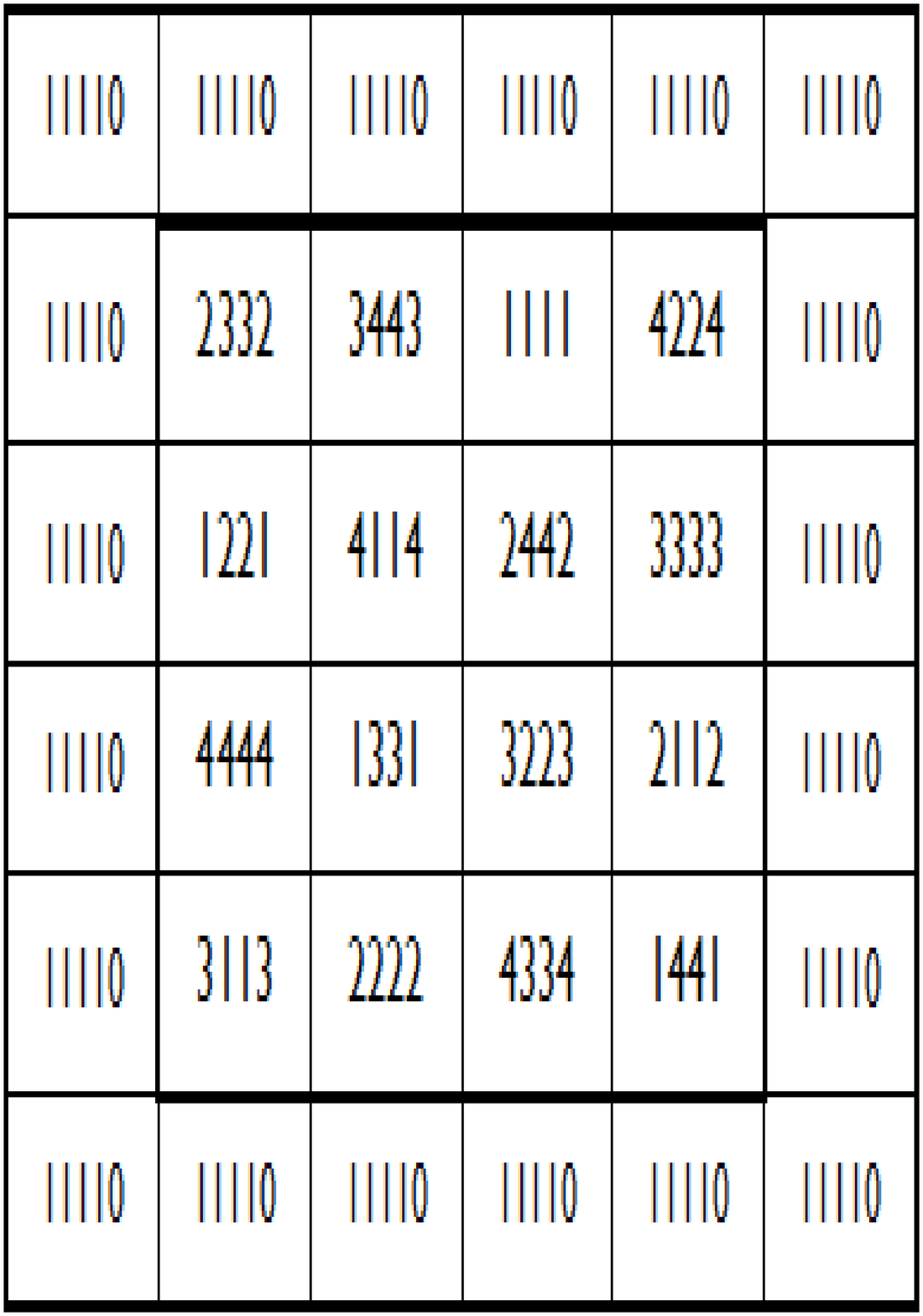}
\end{center}

Writing in digital form we have the following \textbf{upside down magic square} of sum 011110:

\begin{center}
\includegraphics[bb=0mm 0mm 208mm 296mm, width=90.8mm, height=41.3mm, viewport=3mm 4mm 205mm 292mm]{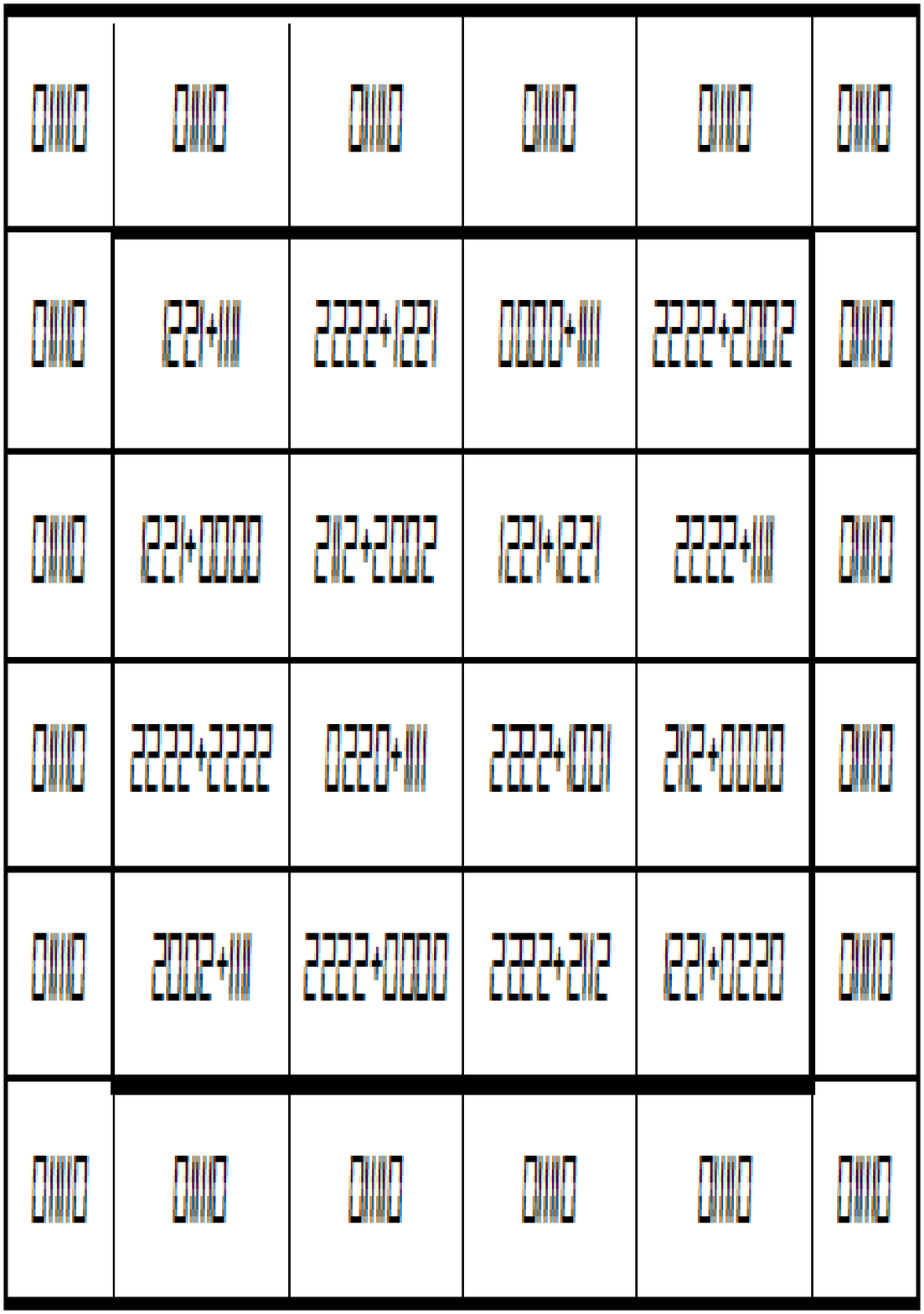}
\end{center}

Here all the numbers are palindromic, except two, the one is 220 and another is the sum 11110. To make them symmetric we have written 0 in the front, i.e, 0220 and 011110.

\subsection{Symmetric and Upside Down Version of ``Lo-Shu'' Magic Square}

We can get ``Lo-Shu'' magic square making an operation $3\times A+B+1$ on the following two \textbf{\textit{non diagonal orthogonal Latin squares}} of order $3 \times 3$

\begin{center}
\includegraphics[bb=0mm 0mm 208mm 296mm, width=32.0mm, height=17.1mm, viewport=3mm 4mm 205mm 292mm]{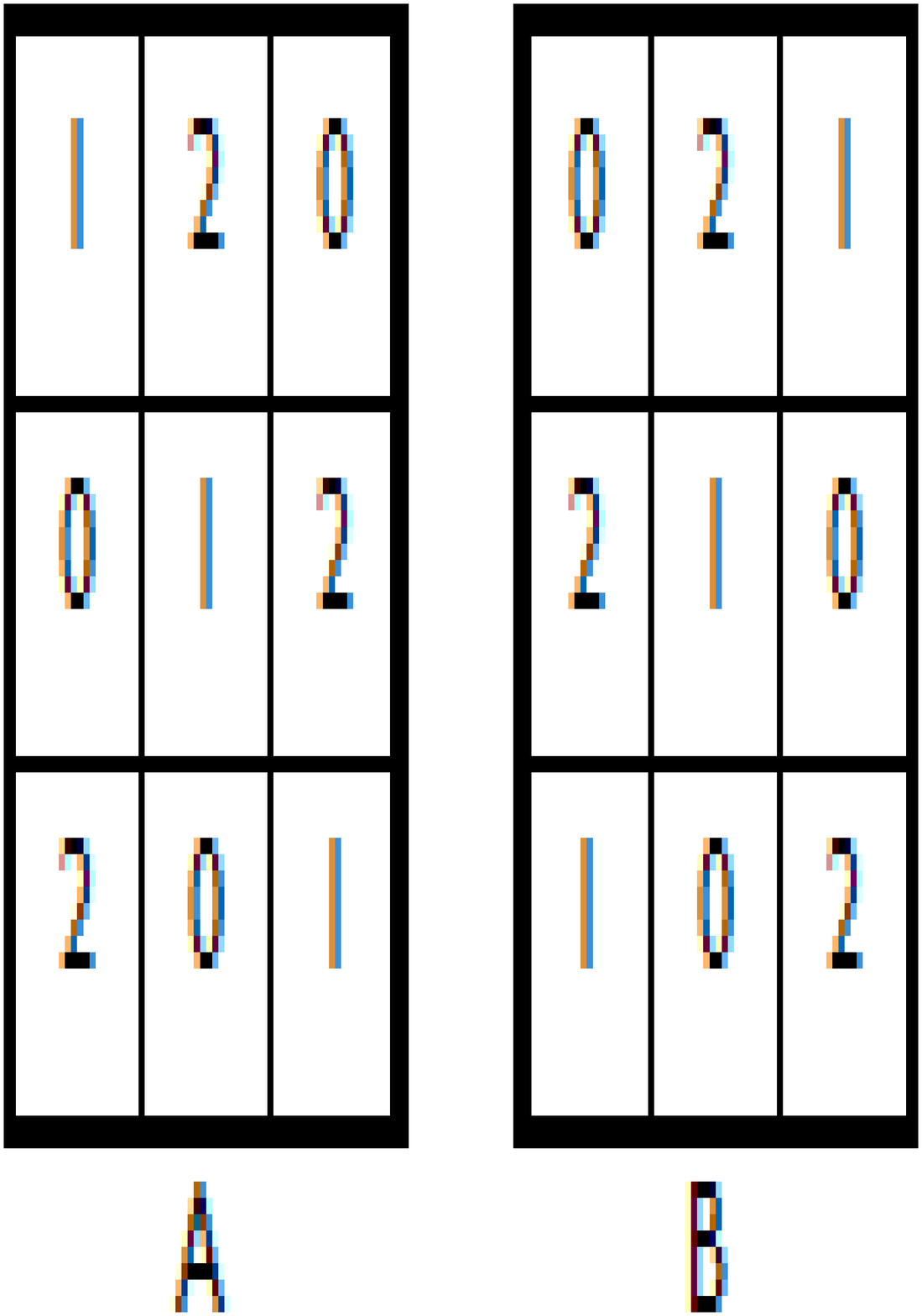}
\end{center}

Again making an operation $10\times A+B$ on the above two Latin squares, we get a 33 sum magic square of order 3. The sum 33 can be written as 11+22. Writing the numbers in the digital form, we have the following \textbf{\textit{equivalent upside down version}} of ``Lo-Shu'' magic square only with the digits 0, 1 and 2.

\begin{center}
\includegraphics[bb=0mm 0mm 208mm 296mm, width=46.8mm, height=32.2mm, viewport=3mm 4mm 205mm 292mm]{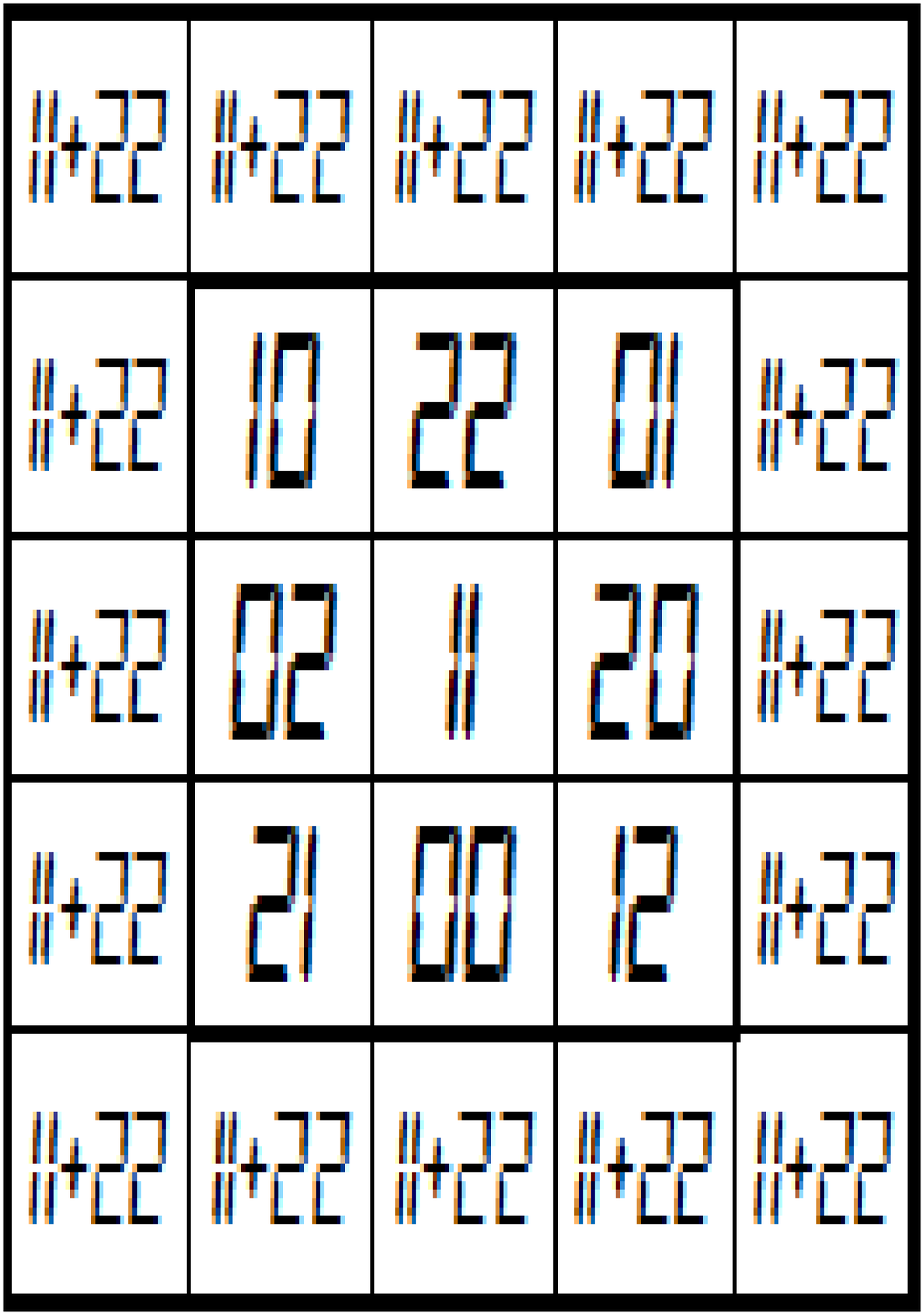}
\end{center}

\subsection{Curiosities }

\begin{itemize}
\item [(i)] The \textit{Khajuraho magic square} was discovered in $10^{th}$  century, and we are in $21^{st}$ century . Both these centuries combined have the digits 0, 1 and 2.

\item [(ii)]  We are in the year 2010. This has only three digits 0, 1 and 2.

\item  [(iii)] The day of submission of this work, i.e., 01.11.2010 also has these three digits 0, 1 and 2.

\item  [(iv)] If we consider the day 01.11.10, i.e, 011110. Still we can have palindromic magic square of sum 011110 only with the digits 0, 1 and 2.

\item  [(v)] The years 2200 BC, 2500 BC, 2800 BC, appearing in the history in different sites of internet regarding the magic square of order $3 \times 3$, have the digits 0-2-5-8. All these four digits written in the digital form \textbf{\textit{\includegraphics[bb=0mm 0mm 208mm 296mm, width=17.4mm, height=3.0mm, viewport=3mm 4mm 205mm 292mm]{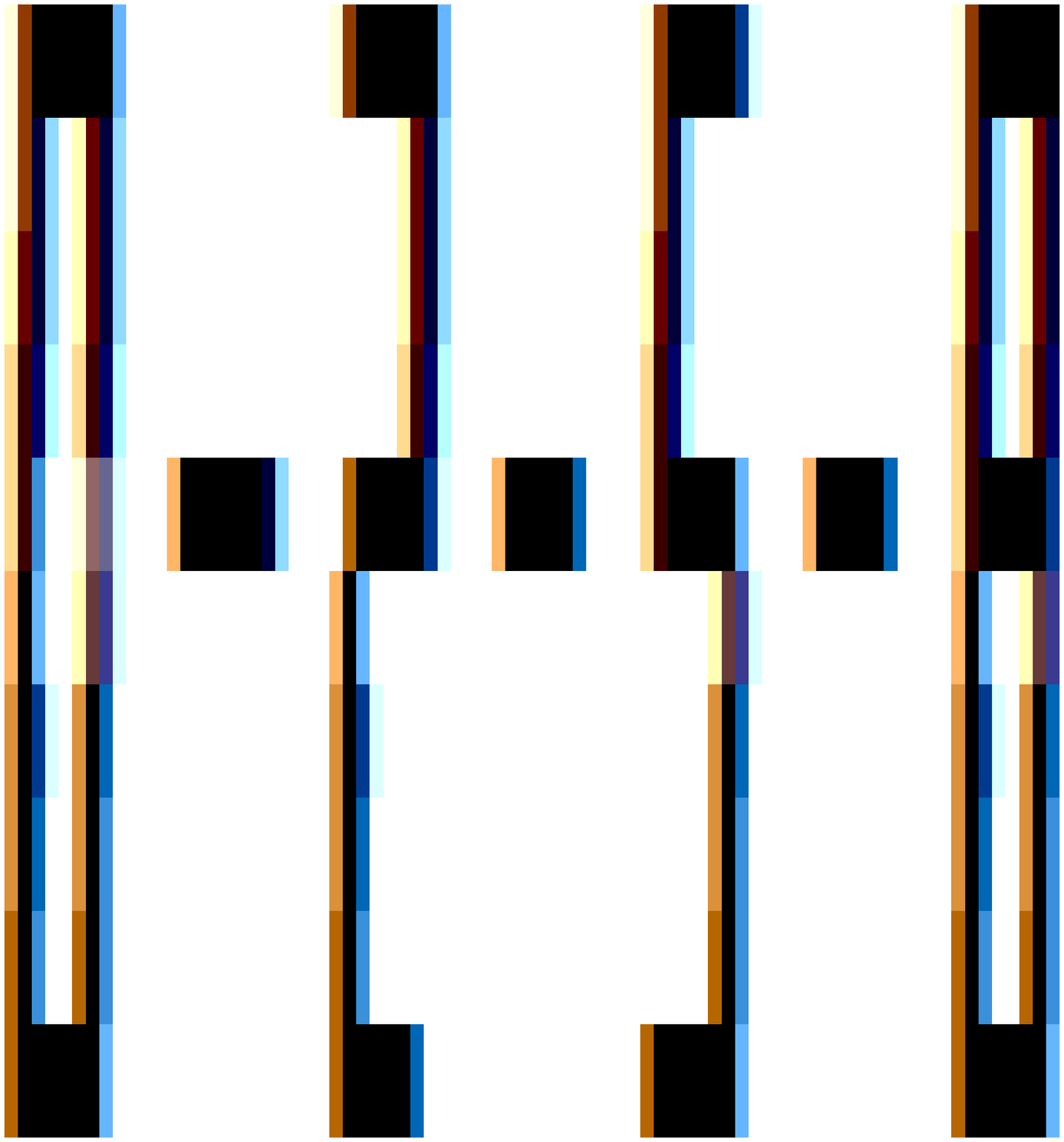}}} are \textit{upside down} and \textit{mirror looking}, where 2 becomes 5 and 5 as 2 when we see them in the mirror.

\item  [(vi)] The magic square $3 \times 3$ is made from the patterns, rather than numbers can be considered as universal, i.e., \textit{upside down} and \textit{mirror looking, }because looking any way it is always the same. See the second plate in section 1.

\item  [(vii)] There are many other days happening during the year 2010, have the digits 0, 1 and 2. Such as 01.01.2010, 02.01.2010, 20.10.2010, 11.11.2010, etc. The day 01.10.2010, we have taken just randomly.
\end{itemize}

In \cite{tan3}, the author gave a bimagic square of order $9 \times 9$ using only the digits 0, 1 and 2, and bimagic square of order $8 \times 8$ using only the digits 0 and 1. For more work on magic square using digital numbers can be seen in Taneja  \cite{tan1, tan2, tan4}. Also the sites \cite{boyer, hei} are good for further studies on magic squares.

\bigskip
The historical figures or plates given above are accessed from the internet, whose links are given in the references \cite{jain, hms, msw}. The historical part is also derived from these sites.

\begin{center}
---------------------------
\end{center}

\end{document}